\documentclass{article}%
\usepackage{amssymb}
\usepackage{amsfonts}
\usepackage{amsmath}
\usepackage{graphicx}%
\providecommand{\U}[1]{\protect\rule{.1in}{.1in}}
\newtheorem{theorem}{Theorem}

\newtheorem{corollary}[theorem]{Corollary}

\newtheorem{definition}[theorem]{Definition}

\newtheorem{lemma}[theorem]{Lemma}

\newtheorem{remark}[theorem]{Remark}

\newenvironment{proof}[1][Proof]{\noindent\textbf{#1.} }{\ \rule{0.5em}{0.5em}}
\begin{document}

\title{Chain recurrence and structure of $\omega$-limit sets of multivalued semiflows}
\author{Oleksiy V. Kapustyan${}^{2}$, Pavlo O. Kasyanov${}^{1}$, Jos\'{e}
Valero${}^{3}$\\$^{1}${\small Institute for Applied System Analysis,}\\{\small National Technical University of Ukraine}\\{\small \textquotedblleft Igor Sikorsky Kyiv Polytechnic
Institute\textquotedblright,}\\{\small Peremogy ave. 37, building 35,}\\{\small 03056, Kyiv, Ukraine}\\{\small kasyanov@i.ua}\\$^{2}${\small Taras Shevchenko National University of Kyiv,}\\{\small Volodymyrska Street 60,}\\{\small 01601, Kyiv, Ukraine}\\{\small alexkap@univ.kiev.ua}\\$^{3}${\small Universidad Miguel Hernandez de Elche,}\\{\small Centro de Investigaci\'{o}n Operativa,}\\{\small Avda. Universidad s/n,}\\{\small 03202-Elche (Alicante), Spain }\\{\small jvalero@umh.es }}
\date{Dedicated to professor Tom\'{a}s Caraballo on the occasion of his 60-th birthday}
\maketitle

\begin{abstract}
We study properties of $\omega$-limit sets of multivalued semiflows like chain
recurrence or the existence of cyclic chains.

First, we prove that under certain conditions the $\omega$-limit set of a
trajectory is chain recurrent, applying this result to an evolution
differential inclusion with upper semicontinous right-hand side.

Second, we give conditions ensuring that the $\omega$-limit set of a
trajectory contains a cyclic chain. Using this result we are able to check
that the $\omega$-limit set of every trajectory of a reaction-diffusion
equation without uniqueness of solutions is an equilibrium.

\end{abstract}

\bigskip

\textbf{AMS Subject Classification (2010):} 35B40, 35B41, 35K55, 37B20, 37B25,
37D15, 58C06

\textbf{Keywords: }multivalued semiflows, $\omega$-limit set, global
attractor, differential inclusions, reaction-diffusion equations, structure,
asymptotic behaviour

\section{Introduction}

The asymptotic behavior of infinite-dimensional dynamical systems without
uniqueness (that is, multivalued dynamical systems) have been intensively
studied during the last three decades (see, among many others, \cite{arv1},
\cite{BabinVishik85}, \cite{ball}, \cite{Ball04}, \cite{Caballero},
\cite{CarLanVal02}, \cite{CakLanVal05}, \cite{CarLanVal2016},
\cite{CarLanVal2019}, \cite{CarMarinRob03}, \cite{CheVi97}, \cite{CostaValero}%
, \cite{CostaVal17}, \cite{Kap1}, \cite{KapPanVal}, \cite{KapVal00},
\cite{KKV}, \cite{KapKasVal15}, \cite{Kas11}, \cite{KloedenVal05}, \cite{li},
\cite{MelVal98}, \cite{SimGen08}, \cite{ZK}).

One important question when studying the asymptotic behavior of solutions of
partial differential equations is to know the internal structure of $\omega
$-limit sets and global attractors, which gives us an insight into the
dynamics of solutions in the long term. While in the single-valued case (for
differential equations with uniqueness of the Cauchy problem) such question
has been widely studied, the multivalued case (for differential equations
without uniqueness of the Cauchy problem) is more difficult to tackle.
Nevertheless, several results in this direction have been published over the
last years (see \cite{arv1}, \cite{ball}, \cite{Caballero}, \cite{CakLanVal05}%
, \cite{CarLanVal2016}, \cite{CarLanVal2019}, \cite{CostaValero},
\cite{CostaVal17}, \cite{KKV}, \cite{KapKasVal15}, \cite{li}).

In this paper we are not interested in studying the structure of the whole
global attractor but the dynamical properties of the $\omega$-limit set of
each particular trajectory of the dynamical system.

First, in Section 2, following the classical theory for single-valued
dynamical systems \cite{Hurley95}, \cite{MishST95}, we introduce the notion of
chain recurrence for multivalued semiflows and prove that under certain
conditions the $\omega$-limit set of an individual trajectory is chain
recurrent. It is worth observing that Conley proved in \cite[Lemma
4.1E]{Conley} that if a point is not chain recurrent for a flow, then a
Lyapunov function exists along its trajectory. Therefore, the property of
being chain recurrent is somehow opposite to the existence of a Lyapunov function.

Second, in Section 3 we apply the abstract theorem of Section 2 to an
evolution inclusion with upper-semicontinuous right-hand side. Moreover, the
converse statement saying that a given compact,connected, quasi-invariant,
chain recurrent set has to be the $\omega$-limit set of a certain differential
inclusion is also established.

Third, in Sections 4-5 we study the internal structure of $\omega$-limit sets,
and in particular the existence of cyclic chains. Such results are very useful
in order to determine whether a trajectory converges towards an equilibrium as
time goes to infinity or not. When a Lyapunov function exists it is possible
to establish, in a similar way as in the single-valued case, that each
trajectory converges to the set of stationary points (see \cite{arv1},
\cite{KKV}). If, moreover, the number of stationary points is finite (or even
infinite but countable), then the $\omega$-limit set of any trajectory is
equal to one stationary point. However, in absence of a Lyapunov function such
results are much harder to prove.

In Section 4 we extend first a classical result \cite{ButlerWaltman} about the
existence of stable and unstable sets for compact, isolated, invariant sets
intersecting with (but not containing) the $\omega$-limit set of one
trajectory. Using it we establish that under certain conditions the $\omega
$-limit set of a trajectory contains a cyclic chain. This theorem generalized
a classical one for semigroups \cite{Thieme}.

In Section 5 we apply these results to a reaction-diffusion equation without
uniqueness of solutions. Although it was proved in \cite{KapKasVal15} that
inside the global attractor the $\omega$-limit set of every trajectory belongs
to the set of stationary points, whether such result is also true or not for
any weak solution of the equation was an open problem so far. Using the
theoretical results of Section 4 we give an answer to this question by proving
that indeed any trajectory converges to an equilibrium if its number is
finite. The idea behind the proof is the following:\ if the $\omega$-limit set
of trajectory was not an equilibrium, then a cyclic chain connecting
equilibria would exist inside the $\omega$-limit set;\ however, as a Lyapunov
function exists in the global attractor, cyclic chains are forbidden.

\section{Chain recurrence for multivalued semiflows}

In this section we will prove that the $\omega$-limit sets of trajectories of
multivalued semiflows are chain recurrent, generalizing in this way the
classical result of Conley for single-valued flows.

Let $X$ be a complete metric space with metric $\rho.$ As usual, the Hausdorff
semidistance from the set $A$ to the set $B$ is given by%
\[
dist\left(  A,B\right)  =\sup_{a\in A}\inf_{b\in B}\ \rho\left(  a,b\right)
.
\]

We consider a set of functions $\mathcal{K}\subset\mathcal{C}(\mathbb{R}%
^{+},X)$ satisfying the following conditions:

\begin{enumerate}
\item[(K1)] for any $x\in X$ there exists $\varphi\in\mathcal{K}$ such that
$\varphi(0)=x$;

\item[(K2)] $\varphi^{\tau}(\cdot)=\varphi(\cdot+\tau)\in\mathcal{K}$, if
$\varphi\in\mathcal{K}$ for any $\tau\geq0$;

\item[(K3)] if $\varphi_{1},\varphi_{2}\in\mathcal{K}$ satisfy $\varphi
_{2}(0)=\varphi_{1}(s)$, $s>0$, then $\varphi$ given by%
\[
\varphi(t)=\left\{
\begin{array}
[c]{ll}%
\varphi_{1}(t), & \mbox{ if }t\leq s\\
\varphi_{2}(t-s), & \mbox{ if }t>s,
\end{array}
\right.
\]
belongs to $\mathcal{K}$;

\item[(K4)] if $\varphi_{n}\in\mathcal{K}$ is a sequence such that
$\varphi_{n}(0)\rightarrow x_{0}$, for some $x_{0}\in X$, then there is a
subsequence and $\varphi_{0}\in\mathcal{K}$ such that $\varphi_{n_{k}%
}(t)\rightarrow\varphi_{0}(t)$ uniformly on compact subsets of $[0,\infty)$.
\end{enumerate}

\begin{remark}
Condition $\left(  K4\right)  $ is stronger than the usual one \cite{ball},
where pointwise convergence is assumed.
\end{remark}

Let $\mathcal{P}(X)$ be the set of all non-empty subsets of $X$. The
multivalued map $G\colon\mathbb{R}^{+}\times X\rightarrow\mathcal{P}(X)$ is
said to be a multivalued semiflow ($m$-semiflow for short) if:

\begin{enumerate}
\item[(i)] $x = G(0,x)$, for all $x\in X$;

\item[(ii)] $G(t+s,x)\subset G(t,G(s,x))$ for all $t,s\geq0$ and $x\in X$.
\end{enumerate}

It is called strict if, additionally, $G(t+s,x)=G(t,G(s,x))$ for all
$t,s\geq0$ and $x\in X$.

We define the multivalued map $G\colon\mathbb{R}^{+}\times X\rightarrow
\mathcal{P}(X)$ associated with the family $\mathcal{K}$ as follows:
\begin{equation}
G(t,x)=\{y\in X:y=\varphi(t)\text{ for some\ }\varphi\in\mathcal{K}\text{ such
that }\varphi(0)=x\}. \label{141014-1640}%
\end{equation}
Conditions $\left(  K1\right)  -\left(  K2\right)  $ imply that $G$ is a
multivalued semiflow;\ if, additionally, $(K3)$ is true, then $G$ is strict
(see e.g. \cite[Lemma 5]{KKV}).

We say that a map $\phi\colon\mathbb{R}\rightarrow X$ is a complete trajectory
of $\mathcal{K}$ if
\[
\phi(\cdot+h)|_{[0,\infty)}\in\mathcal{K},\text{ for any }h\in\mathbb{R}.
\]

We recall several definitions of invariance for a set $A\subset X.$ $A$ is
said to be positively invariant if $G(t,A)\subset A,$ for all $t\geq0,$ and
negatively invariant if $A\subset G(t,A)$ for all $t\geq0$. It is invariant if
it is both positively and negatively invariant, that is, $G(t,A)=A$ for all
$t\geq0$. $A$ is quasi-invariant (or weakly invariant as well) if for all
$x\in A$ there is at least one complete trajectory $\phi$ of $\mathcal{K}$
such that $\phi(t)\in A$, for all $t\in\mathbb{R}$.

It is obvious that any quasi-invariant set is negatively invartiant. It is
also well known that under conditions $\left(  K1\right)  -\left(  K4\right)
$ any compact invariant set is quasi-invariant \cite[Corollary 7]{CostaValero}.

For any trajectory $\varphi\in\mathcal{K}$ we define its $\omega$-limit set by%
\begin{align*}
\omega\left(  \varphi\right)   &  =\underset{{\tau\geq0}}{\bigcap}%
\,\overline{\underset{{t\geq\tau}}{\bigcup}\varphi(t)}\\
&  =\{y\in X:\text{there exists a sequence }t_{n}\rightarrow+\infty\text{ such
that }\varphi(t_{n})\rightarrow y\}.
\end{align*}
The positive orbit of $\varphi\in\mathcal{K}$ is the set $\gamma^{+}%
(\varphi)=\cup_{t\geq0}\varphi\left(  t\right)  $.

If $\phi$ is a complete trajectory of $\mathcal{K}$, the $\alpha$-limit set is
defined by%
\begin{align*}
\alpha\left(  \phi\right)   &  =\underset{{\tau\leq0}}{\bigcap}\,\overline
{\underset{{t\leq\tau}}{\bigcup}\phi(t)}\\
&  =\{y\in X:\text{there exists a sequence }t_{n}\rightarrow-\infty\text{ such
that }\phi(t_{n})\rightarrow y\}.
\end{align*}
The negative orbit of $\phi$ is the set $\gamma^{-}(\phi)=\cup_{t\leq0}%
\phi\left(  t\right)  $.

The following lemma can be proved in the same way as in Lemma 3.4 and
Proposition 4.1 in \cite{ball}.

\begin{lemma}
\label{PropOmega}Let $\left(  K1\right)  -\left(  K4\right)  $ be satisfied.
If the closure of the positive orbit of $\varphi\in\mathcal{K}$ is compact,
then $\omega\left(  \varphi\right)  $ is non-empty, compact, connected,
quasi-invariant and
\[
\lim_{t\rightarrow+\infty}dist\left(  \varphi\left(  t\right)  ,\omega\left(
\varphi\right)  \right)  =0.
\]
If $\phi$ is a complete trajectory such that the closure of the negative orbit
is compact, then $\alpha(\phi)$ is non-empty, compact, connected,
quasi-invariant and
\[
\lim_{t\rightarrow-\infty}dist\left(  \phi\left(  t\right)  ,\alpha\left(
\phi\right)  \right)  =0.
\]

\end{lemma}

\begin{definition}
\label{defi:1} Let $B\subset H$ be a quasi-invariant set with respect to
m-semiflow $G$, and let $y,z\in B$. For $\varepsilon>0$, $t>0$ an
$(\varepsilon,t)$-chain from $y$ to $z$ is a sequence $\{y=y_{1}%
,y_{2},...,y_{n+1}=z\}\subset B$, $\{t_{1},t_{2},...,t_{n}\}\subset\lbrack
t,+\infty)$ such that
\begin{equation}
\mathrm{dist}(y_{i+1},G(t_{i},y_{i}))<\varepsilon,\,\,\,i=\overline{1,n}.
\label{dcr}%
\end{equation}
A point $y\in A$ is called chain recurrent with respect to $G$ if for every
$\varepsilon>0$, $t>0$ there exists an $(\varepsilon,t)$-chain from $y$ to
$y$. The set $B$ is said to be chain recurrent with respect to $G$ if every
point of $B$ is chain recurrent with respect to $G$.
\end{definition}

\begin{remark}
If $G$ is single-valued, then quasi-invariance implies positively invariance
and, as a consequence, Definition \ref{defi:1} coincides with the classical
definition of chain recurrence \cite{MishST95}.
\end{remark}

Let us consider the following conditions:

\begin{itemize}
\item[$\left(  N1\right)  $] There exist a sequence of sets of functions
$\mathcal{K}_{N}\subset\mathcal{C}(\mathbb{R}^{+},X)$ satisfying $\left(
K1\right)  -\left(  K4\right)  $ such that:
\begin{equation}
\ G(t,y)\subset G_{N}(t,y),\ \text{for all }t\geq0,\ N\geq1,\ y\in X,
\label{vv_eq_12}%
\end{equation}%
\begin{equation}
\mathrm{dist}_{H}\left(  G_{N}(t,y_{1}),G_{N}(t,y_{2})\right)  \leq e^{c_{N}%
t}\Vert y_{1}-y_{2}\Vert,\text{ for all }y_{1},y_{2}\in X,\ t\geq0,\ N\geq1,
\label{vv_eq_13}%
\end{equation}
where $G_{N}:\mathbb{R}^{+}\times X\rightarrow\mathcal{P}(X)$ are the strict
m-semiflows corresponding to $\mathcal{K}_{N}$ and $\mathrm{dist}_{H}\left(
A,B\right)  =\max\{\mathrm{dist}\left(  A,B\right)  ,\mathrm{dist}(B,A)\}$ is
the Hausdorff distance.

\item[$\left(  N2\right)  $] If $y_{N}\in\mathcal{K}_{N}$, $y_{N}%
(0)\rightarrow y_{0}$, then up to subsequence
\begin{equation}
y_{N}\rightarrow y\ \text{in}\ C([0,T];X),\text{ for any }T>0,
\label{vv_eq_6++}%
\end{equation}
where $y\in\mathcal{K}$, $y(0)=y_{0}$.
\end{itemize}

\begin{theorem}
\label{th:1} Assume that conditions $\left(  K1\right)  -\left(  K4\right)  $,
$\left(  N1\right)  -\left(  N2\right)  $ hold. If the closure of the positive
orbit of $\varphi\in\mathcal{K}$ is compact, then the set $\omega(\varphi)$ is
chain recurrent with respect to the strict m-semiflow $G$.
\end{theorem}

\begin{proof}
It is well known that $\varphi(t+s)\in G(t,\varphi(s))$ \cite{KapKasVal15}.
Let us fix $T$ and let us put
\[
\gamma^{T}(\varphi)=\bigcup\limits_{t\geq T}\varphi(t).
\]
We take arbitrary $\varepsilon>0$, $y\in\omega(\varphi)$, $t_{0}>0$. Because
of the equality
\[
y=\lim\varphi(s_{n})\ \text{as}\ s_{n}\rightarrow\infty,
\]
and $\left(  K4\right)  $ we can choose $n$ such that $s_{n}>T$ and%
\[
\ \mathrm{dist}\left(  \varphi(s_{n}+t),G(t,y)\right)  <\varepsilon
\,,\,\text{for all}\ t\in\lbrack0,t_{0}].
\]
Let us put $y_{1}=y$, $y_{2}=\varphi(s_{n}+t_{0})$, $t_{1}=t_{0}$. Then
\[
\ \mathrm{dist}(y_{2},G(t_{1},y_{1}))=\mathrm{dist}(\varphi(s_{n}%
+t_{0}),G(t_{0},y))<\varepsilon.
\]
Choose $m\geq n$ such that
\[
s_{m}>s_{n}+2t_{0},\,\,\,\Vert\varphi(s_{m})-y\Vert<\varepsilon.
\]
Let $k\geq1$ be such that
\[
s_{m}-s_{n}-t_{0}=kt_{0}+r,\,\,\,\ r\in\lbrack0,t_{0}).
\]
Let us put
\[
y_{3}=\varphi(s_{n}+2t_{0}),\ldots,\,\ y_{k+1}=\varphi(s_{n}+kt_{0}%
),\,\ y_{k+2}=y,
\]%
\[
t_{1}=\ldots=t_{k}=t_{0},\,\,\ t_{k+1}=t_{0}+r.
\]
Then
\[
\mathrm{dist}\left(  y_{i+1},G(t_{i},y_{i})\right)  =\mathrm{dist}\left(
\varphi(s_{n}+it_{0}),G(t_{0},\varphi(s_{n}+(i-1)t_{0}))\right)  =0,\text{ for
any }i\in\lbrack2,k],
\]%
\[
\mathrm{dist}(y_{k+2},G(t_{k+1},y_{k+1}))=\mathrm{dist}(y,G(t_{0}%
+r,\varphi(s_{n}+kt_{0})))\leq\Vert y-\varphi(s_{m})\Vert<\varepsilon.
\]
In other words, for any $y\in\omega(\varphi),$ $\varepsilon>0$, $t_{0}>0$
there exists an $(\varepsilon,t_{0})$-chain $\{y=y_{1},\ldots,y_{l+1}=y\}$,
$\{t_{1},\ldots,t_{l}\}$ such that $y_{i}\in\gamma^{T}(\varphi)$,
$i=\overline{2,l}$, $t_{l}\in\lbrack t_{0},2t_{0})$.

From (\ref{vv_eq_12}) the same is true for $G_{N}$ for every $N\geq1$. So,
putting $\varepsilon=\frac{1}{n}$, $T=n$ we obtain that for every $n\geq1$
there exist $\{y_{i}^{n}\}_{i=1}^{l_{n}+1}\subset\gamma^{n}(\varphi)$,
$y_{1}^{n}=y=y_{l_{n}+1}^{n}$, $\{t_{i}^{n}\}_{i=1}^{l_{n}}$, $t_{i}^{n}%
=t_{0}$, $i=1,2,\ldots,l_{n}-1$, $t_{l_{n}}^{n}\in\lbrack t_{0},2t_{0})$ such
that for any $N\geq1,$
\[
\mathrm{dist}(y_{i+1}^{n},G_{N}(t_{i}^{n},y_{i}^{n}))<\frac{1}{n}%
,\,\,\,i={1,2,\ldots,l_{n}-1}.
\]
Let us denote
\[
C^{n}=\bigcup\limits_{i=1}^{l_{n}}y_{i}^{n}.
\]
Then $C^{n}$ is compact, $y\in C^{n}$, $C^{n}\subset\gamma(\varphi
)=\bigcup\limits_{t\geq0}\varphi(t),\ $for all $n\geq1$. Due to the
compactness of $\overline{\gamma(\varphi)}$ up to a subsequence
\[
\mathrm{dist}_{H}(C^{n},C)\rightarrow0,\,\,\ n\rightarrow\infty,
\]
where $y\in C$, $C\subset\omega(\varphi)$. For every $N\geq1$ we choose
$n_{0}$ such that for all $n\geq n_{0},$
\[
\frac{1}{n}<\frac{\varepsilon}{3},\ \alpha_{n}:=\mathrm{dist}_{H}%
(C^{n},C)<\frac{\varepsilon}{3}\ \ \text{and}\ \ e^{2c_{N}t_{0}}\alpha
_{n}<\frac{\varepsilon}{3}.
\]
Also let us put $z_{1}=y$, $t_{1}=t_{0}$, $z_{2}^{n}\in C$ such that
\[
\mathrm{dist}(y_{2}^{n},C)=\mathrm{dist}(y_{2}^{n},z_{2}^{n}).
\]
Then
\[
\mathrm{dist}(z_{2}^{n},G_{N}(t_{1},z_{1}))\leq\mathrm{dist}(z_{2}^{n}%
,y_{2}^{n})+\mathrm{dist}(y_{2}^{n},G_{N}(t_{0},y))<\frac{2\varepsilon}%
{3}<\varepsilon.
\]
Let us put $t_{2}=t_{0}$, $z_{3}^{n}\in C$ such that $dist(y_{3}%
^{n},C)=dist(y_{3}^{n},z_{3}^{n}).$ Then
\[
\mathrm{dist}(z_{3}^{n},G_{N}(t_{2},z_{2}^{n}))\leq\mathrm{dist}(z_{3}%
^{n},y_{3}^{n})+\mathrm{dist}(y_{3}^{n},G_{N}(t_{0},y_{2}^{n}))+
\]%
\[
+\mathrm{dist}(G_{N}(t_{0},y_{2}^{n}),G_{N}(t_{0},z_{2}^{n}))<\frac
{\varepsilon}{3}+\frac{1}{n}+\frac{\varepsilon}{3}<\varepsilon.
\]
Repeating this argument we obtain $z_{i}^{n}\in C$, $i\in\{1,...,l_{n}\}$,
$t_{i}=t_{0}$, for $i=1,...,l_{n}-1,$ and%
\[
\mathrm{dist}(z_{i+1}^{n},G_{N}(t_{i},z_{i}^{n}))<\varepsilon.
\]
Let us put $z_{l_{n}+1}=y$, $t_{l_{n}}^{n}\in\lbrack t_{0},2t_{0})$. Then
\[
\mathrm{dist}(z_{l_{n}+1},G_{N}(t_{l_{n}}^{n},z_{l_{n}})\leq\mathrm{dist}%
(y,G_{N}(t_{l_{n}}^{n},y_{l_{n}}^{n}))+\mathrm{dist}(G_{N}(t_{l_{n}}%
^{n},y_{l_{n}}^{n}),G_{N}(t_{l_{n}}^{n},z_{l_{n}}))<\varepsilon.
\]
Due to (\ref{vv_eq_6++}) and the compactness of the set $C$ we can choose a
number $N$ such that
\[
\sup\limits_{t\in\lbrack t_{0},2t_{0}]}\sup\limits_{z\in C}\mathrm{dist}%
(G_{N}(t,z),G(t,z))<\varepsilon
\]
Then
\[
\mathrm{dist}(z_{i+1}^{n},G(t_{i},z_{i}^{n}))\leq\mathrm{dist}(z_{i+1}%
^{n},G_{N}(t_{i},z_{i}^{n}))+\varepsilon<2\varepsilon
\]
and the theorem is proved.
\end{proof}

\bigskip

As a consequence of (\ref{dcr}) we have the following result, which can be
proved by repeating without any changes the arguments of Lemmas 1.4, 3.3 from
\cite{MishST95}.

\begin{corollary}
\label{cor:1} Assume that conditions $\left(  K1\right)  -\left(  K4\right)  $
hold. Let $B$ be a compact connected chain recurrent set with respect to the
m-semiflow $G$. Then for any $T>0$,$\ \varepsilon>0$, $y_{0}\in B$ there exist
sequences $\{y_{i}\}_{i\geq1}\subset B$, $\{t_{i}\}_{i\geq0}\subset\lbrack
T,+\infty)$ such that
\begin{equation}
\mathrm{dist}(y_{i+1},G(t_{i},y_{i}))<\varepsilon\ \text{for all }i\geq0,
\label{m2}%
\end{equation}%
\begin{equation}
\mathrm{dist}(y_{i+1},G(t_{i},y_{i}))\rightarrow0,\ \text{as }i\rightarrow
\infty, \label{m3}%
\end{equation}%
\begin{equation}
B=\overline{\{y_{i}\}}_{i\geq n}\ \text{for all }\ n\geq1. \label{m4}%
\end{equation}

Let us denote
\[
D(y_{0})=\{y(\cdot)\in\mathcal{K}\ |\ y(0)=y_{0}\}.
\]
Then for any $i\geq0$ there exists $\varphi_{i}\in D(y_{i})$ such that
\[
\mathrm{dist}(y_{i+1},G(t_{i},y_{i}))\leq\Vert y_{i+1}-\varphi_{i}(t_{i}%
)\Vert=:\varepsilon_{i}<\varepsilon,\ \varepsilon_{i}\searrow0,\ i\rightarrow
\infty.
\]

\end{corollary}

Corollary~\ref{cor:1}, and mappings $\{\varphi_{i}\}_{i\geq0}$ allow us to
construct a pseudo-trajectory, the $\omega$-limit set of which contains $B$.

\begin{corollary}
\label{cor:2}Assume the conditions of Corollary~\ref{cor:1}. Then for any
$\varepsilon>0,\ T>0,\ y_{0}\in B$ there exists a mapping (an $(\varepsilon
,T)$-pseudo-trajectory) $\varphi^{\ast}(\cdot)$ starting at $y_{0}$ such that
\begin{equation}
B\subset\omega(\varphi^{\ast}). \label{m7}%
\end{equation}
This mapping is defined by the following formula
\begin{equation}
\varphi^{\ast}(t)=\left\{
\begin{array}
[c]{l}%
\varphi_{0}(t)\in G(t,y_{0}),\ t\in\lbrack0,t_{0}),\\
\varphi_{i}(t-s_{i-1})\in G(t-s_{i-1},y_{i}),\ s_{i-1}\leq t<s_{i}.
\end{array}
\right.  \label{m5}%
\end{equation}
where $s_{i}=\sum\limits_{k=0}^{i}t_{k}$ and
\begin{equation}
\varepsilon_{i}=\Vert\varphi^{\ast}(s_{i})-\varphi^{\ast}(s_{i}-0)\Vert
\rightarrow0\text{.} \label{m6}%
\end{equation}

\end{corollary}

\section{Chain recurrence for differential inclusions\label{Incl}}

Let $V\subset H\subset V^{\prime}$ be a Gelfand triple with compact dense
embeddings, $\Vert\cdot\Vert$ and $(\cdot,\cdot)$ be the norm and the scalar
product in $H$, $\langle\cdot,\cdot\rangle$ be pairing between $V$ and
$V^{\prime}$. We are interested in the limit behavior of trajectories of the
following evolution inclusion
\begin{equation}
\left\{
\begin{array}
[c]{lll}%
\dfrac{dy}{dt}+Ay\in F(y), &  & \\
y|_{t=0}=y_{0}\in H, &  &
\end{array}
\right.  \label{vv_eq_1}%
\end{equation}
where $A:V\rightarrow V^{\prime}$ is a linear, continuous, self-adjoint
operator satisfying
\begin{equation}
\lambda_{1}\Vert u\Vert_{V}^{2}\leq\langle Au,u\rangle\leq\lambda_{2}\Vert
u\Vert_{V}^{2}\text{, }0<\lambda_{1}<\lambda_{2}. \label{vv_eq_2}%
\end{equation}
and the multi-valued term $F:H\rightarrow\mathcal{P}(H)$ satisfies the
following assumptions:
\begin{align}
&  F\text{ has closed, convex, bounded values,}\label{vv_eq_3}\\
&  F\text{ is w-upper semicontinuous,}\nonumber\\
\left\Vert F(y)\right\Vert _{+}  &  :=\sup_{a\in F(y)}\left\Vert a\right\Vert
\leq C_{1}+C_{2}\left\Vert y\right\Vert \text{, for all }y\in H,\nonumber
\end{align}
for some constants $C_{1},C_{2}>0.$

We note that the continuous embedding $V\subset H$ and (\ref{vv_eq_2}) imply
that%
\[
\langle Au,u\rangle\geq\gamma\Vert u\Vert^{2}\text{, for some }\gamma
>0\text{.}%
\]

The assumptions on $A$ imply that it is in fact a densely defined maximal
monotone operator in $H$, so it is well known \cite[p.60]{Barbu} that $A$ is
equal to the subdifferential $\partial\psi$ of the proper, convex lower
semicontinuous function%
\[
\psi\left(  u\right)  =\left\{
\begin{array}
[c]{c}%
\frac{1}{2}\left\Vert A^{\frac{1}{2}}u\right\Vert \text{, if }u\in
D(A^{\frac{1}{2}}),\\
+\infty\text{, otherwise.}%
\end{array}
\right.
\]
As $\left\Vert A^{\frac{1}{2}}u\right\Vert $ is an equivalent norm in $V$, the
compact embedding $V\subset H$ implies also that the sets
\begin{equation}
M_{R}=\{u\in H\mid\left\Vert u\right\Vert \leq R,\ \psi\left(  u\right)  \leq
R\} \label{Level}%
\end{equation}
are compact for any $R>0$.

We recall that $F$ is called w-upper semicontinuous at $y_{0}$ if for any
$\varepsilon>0$ there exists $\delta>0$ such that $F\left(  y\right)  \subset
O_{\varepsilon}(F(y_{0}))$ as soon as $\left\Vert y-y_{0}\right\Vert <\delta,$
where for a set $A\subset H$ we define its $\varepsilon$-neighborhood by
\[
O_{\varepsilon}(A)=\{z\in H:dist\left(  z,A\right)  <\varepsilon\}.
\]
$F$ is w-upper semicontinuous if it is w-upper semicontinuous at any $y_{0}\in
H.$ $F$ is called upper semicontinuous if in the above definition we replace
the $\varepsilon$-neighborhood $O_{\varepsilon}(F(y_{0}))$ by an arbitrary
neighborhood $O(F(y_{0}))$. It is obvious that any upper semicontinuous map is
w-upper semicontinuous, the converse being true as well when $F$ possesses
compact values.

It is known \cite{Den93} that under conditions (\ref{vv_eq_2}), (\ref{vv_eq_3}%
) for every $y_{0}\in H$ problem (\ref{vv_eq_1}) has at least one (mild)
solution $y=y(t)$, $t\geq0$, i.e., for any $T>0$ $y\in L^{2}(0,T;V)$,
$\dfrac{dy}{dt}\in L^{2}(0,T;V^{\prime})$ and there exists $f\in L^{2}(0,T;H)$
such that
\begin{equation}
\left\{
\begin{array}
[c]{ll}%
\dfrac{dy}{dt}+Ay=f(t),\,\,\,f(t)\in F(y(t))\text{ for a.a. }t\in\left(
0,T\right)  , & \\
y|_{t=0}=y_{0}. &
\end{array}
\right.  \label{vv_eq_4}%
\end{equation}

Let $\mathcal{K}\subset C([0,+\infty);H)$ be the collection of all mild
solutions of (\ref{vv_eq_1}). It is known \cite{KapVal00} that properties
$\left(  K1\right)  -\left(  K4\right)  $ are satisfied. Then we define the
strict m-semiflow $G:\mathbb{R}^{+}\times H\rightarrow\mathcal{P}(H)$ by
(\ref{141014-1640}).

We will use the following condition: there exists a sequence of mappings
$\{F_{N}:H\rightarrow\mathcal{P}(H)\}$ such that $F_{N}$ have closed, convex,
bounded values and
\begin{equation}
\,\ F(y)\subset F_{N}(y),\,\,\text{for all }y\in H,\ N\geq1, \label{vv_eq_8}%
\end{equation}%
\begin{equation}
\,\,\mathrm{dist}(F_{N}(y),F(y))\rightarrow0,\,\,\text{as}\,N\rightarrow
\infty,\text{ for all }y\in H, \label{vv_eq_9}%
\end{equation}%
\begin{equation}
\,\mathrm{dist}_{H}(F_{N}(y_{1}),F_{N}(y_{2}))\leq c_{N}\Vert y_{1}-y_{2}%
\Vert\text{, for all }y_{1},y_{2}\in H, \label{vv_eq_10}%
\end{equation}
for some $c_{N}>0$. Let us consider the evolution inclusion
\begin{equation}
\left\{
\begin{array}
[c]{ll}%
\dfrac{dy}{dt}+Ay\in F_{N}(y), & \\
y|_{t=0}=y_{0}\in H. &
\end{array}
\right.  \label{vv_eq_11}%
\end{equation}
All mild solutions $K_{N}$ of (\ref{vv_eq_11}) generate the strict m-semiflows
$G_{N}:\mathbb{R}^{+}\times H\rightarrow\mathcal{P}(H)$ (see \cite{MelVal98})
and (\ref{vv_eq_12})-(\ref{vv_eq_6++}) hold true. Property (\ref{vv_eq_12}) is
obvious, whereas (\ref{vv_eq_13}) was proved in \cite[Lemma 8]{MelVal98}.
Finally, for (\ref{vv_eq_6++}) see the proof of Theorem 3.1 in \cite{KapVal00}.

Therefore, conditions $\left(  K1\right)  -\left(  K4\right)  ,\ \left(
N1\right)  -\left(  N2\right)  $ are satisfied and from Theorem \ref{th:1}\ we
obtain the following result.

\begin{theorem}
\label{ChainIncl}Asume that (\ref{vv_eq_8})-(\ref{vv_eq_10}) are satisfied. If
the closure of the positive orbit of $\varphi\in\mathcal{K}$ is compact, then
the set $\omega(\varphi)$ is chain recurrent with respect to the strict
m-semiflow $G$.
\end{theorem}

Let us give an additional condition ensuring that the closure of every
positive orbit of $\varphi\in\mathcal{K}$ is compact. For this aim we need to
apply a result on existence of global attractors given in \cite{KapVal00}. We
recall that the set $\Theta$ is a global attractor for $G$ if $\Theta\subset
G(t,\Theta)$ for all $t\geq0$ (negatively semi-invariance),
$dist(G(t,B),\Theta)\rightarrow0,$ as $t\rightarrow+\infty,$ for any bounded
ser $B$ (attraction property) and it is minimal (that is, it is contained in
any closed attracting set). It is called invariant if, moreover,
$\Theta=G(t,\Theta)$ for all $t\geq0.$

If we suppose also that for some $\delta>0$,
\begin{equation}
\left(  z,u\right)  \leq\left(  \lambda_{1}-\delta\right)  \left\Vert
u\right\Vert ^{2}\text{, for any }u\in H\text{, }z\in F(u), \label{Diss}%
\end{equation}
then $G$ has the invariant compact global attractor $\Theta$ \cite{KapVal00}.
In particular, every positive orbit $\gamma^{+}(\varphi)$, $\varphi
\in\mathcal{K},$ has compact closure in $H$ and by Lemma \ref{PropOmega}\ its
$\omega$-limit set $\omega(\varphi)$ is nonempty, compact, quasi-invariant,
connected and $\omega(\varphi)\subset\Theta$.

\begin{corollary}
\label{ChainInclCor}Asume that (\ref{vv_eq_8})-(\ref{vv_eq_10}) are satisfied.
If (\ref{Diss}) holds, then for any $\varphi\in\mathcal{K}$ the set
$\omega(\varphi)$ is chain recurrent with respect to the strict m-semiflow $G$.
\end{corollary}

Now let us prove in some sense the converse statement to the previous theorem:
under some additional restrictions on problem (\ref{vv_eq_1}) every chain
recurrent quasi-invariant compact connected set of the m-semiflow $G$ is the
$\omega$-limit set of some trajectory of a perturbed inclusion.

\begin{theorem}
\label{th:2}Asume that (\ref{vv_eq_8})-(\ref{vv_eq_10}) are satisfied. Let $B$
be a compact, connected, chain recurrent, quasi-invariant set with respect to
the m-semiflow $G$. Assume that $A\in\mathcal{L}(H)$. Then there exist maps
$F^{N}:H\rightarrow\mathcal{P}(H)$ satisfying (\ref{vv_eq_3}) and
\[
\ \mathrm{dist}_{H}(F^{N}(y),F(y))\rightarrow0,\ \text{as}\ N\rightarrow
\infty,\text{ for all }y\in H,
\]
such that for any $y_{0}\in B,$ $N\geq1$ there is a solution $y_{N}(\cdot)$ of
the problem
\begin{equation}
\left\{
\begin{array}
[c]{ll}%
\dfrac{dy}{dt}+Ay\in F^{N}(y), & \\
y(0)=y_{0}, &
\end{array}
\right.  \label{vv_eq_14}%
\end{equation}
for which $B\subset\omega(y_{N})$.
\end{theorem}

\begin{proof}
Let us fix $N\geq1$, $\varphi(\cdot)\in K$, $\varphi_{0}\in H$ and
$\overline{t}>1$. Denote by $\left\Vert A\right\Vert _{\mathcal{L}}$ the norm
of $A$ in $\mathcal{L}(H)$.

Let us show that for any $\varepsilon\in(0,\frac{1}{N(c_{N}+\Vert
A\Vert_{\mathcal{L}}+1)}),\ b\in O_{\varepsilon}(\varphi(\overline{t}))$ there
exists a solution $y_{\varepsilon}(\cdot)$ of the problem
\begin{equation}
\left\{
\begin{array}
[c]{ll}%
\dfrac{dy}{dt}+Ay\in F^{N}(y),\ t\in(0,\bar{t}) & \\
y(0)=\varphi_{0}, &
\end{array}
\right.  \label{vv_eq_14+}%
\end{equation}
with $y_{\varepsilon}(\overline{t})=b$, where
\[
F^{N}(y)=\overline{O_{\frac{1}{N}}(F_{N}(y))}.
\]
Indeed, $\varphi(\cdot)\in K_{N}$. We put $\xi=b-\varphi(\overline{t})$ and
consider
\[
y_{\varepsilon}(t)={\varphi}(t)+\xi\frac{t}{\overline{t}},\ t\in
\lbrack0,\overline{t}].
\]
Then $y_{\varepsilon}(0)=\varphi_{0},\ y_{\varepsilon}(\overline{t})=b$
\[
\dot{y}_{\varepsilon}(t)+Ay_{\varepsilon}\in A\xi\frac{t}{\overline{t}}%
+\xi\frac{1}{\overline{t}}+F_{N}(y_{\varepsilon}(t)-\xi\frac{t}{\overline{t}%
}).
\]
Due to (\ref{vv_eq_10}) we deduce that for a.a. $t\in(0,\bar{t})$
\[
\dot{y}_{\varepsilon}(t)+Ay_{\varepsilon}\in O_{(c_{N}+\Vert A\Vert
_{\mathcal{L}}+1)\varepsilon}(F_{N}(y_{\varepsilon}(t))\subset F^{N}%
(y_{\varepsilon}(t)).
\]
Let us denote this solution by $y_{\varepsilon}(t,\varphi_{0})$.

We observe that the maps $F^{N}$ satisfy the conditions in (\ref{vv_eq_3}). It
is clear that $F^{N}$ has closed, convex, bounded values. Also, if we prove
that%
\begin{equation}
dist_{H}(F^{N}(x),F^{N}(y))\leq c_{N}\left\Vert x-y\right\Vert ,\text{ for all
}x,y\in H\text{,} \label{Lip}%
\end{equation}
then, in particular, $F^{N}$ are w-upper semicontinous. Indeed, first by
(\ref{vv_eq_10})\ we have%
\[
dist(F_{N}(x),F^{N}(y))\leq dist(F_{N}(x),F_{N}(y))\leq c_{N}\left\Vert
x-y\right\Vert .
\]
If $z\in F^{N}(x)\backslash F_{N}(x)$, then we choose $v\in F_{N}(x)$
satisfying $z=v+w$, where $\left\Vert w\right\Vert <\frac{1}{N}$. Take an
arbitrary $\delta>0$. There is some $u\in F_{N}(y)$ such that%
\[
\left\Vert v-u\right\Vert \leq dist\left(  v,F_{N}(y))\right)  +\delta\leq
c_{N}\left\Vert x-y\right\Vert +\delta.
\]
Then $\widetilde{u}=u+w\in F^{N}(y)$ and
\[
\left\Vert z-\widetilde{u}\right\Vert =\left\Vert v-u\right\Vert \leq
c_{N}\left\Vert x-y\right\Vert +\delta.
\]
Since $\delta$ and $z$ are arbitrary, we obtain that%
\[
dist\left(  F^{N}(x),F^{N}(y)\right)  \leq c_{N}\left\Vert x-y\right\Vert .
\]
Arguing the other way round we obtain%
\[
dist\left(  F^{N}(y),F^{N}(x)\right)  \leq c_{N}\left\Vert x-y\right\Vert
\]
as well and so (\ref{Lip}). Finally, from (\ref{Lip}) the existence of
constants $C_{1}^{N},C_{2}^{N}>0$ such that
\[
\left\Vert F^{N}(y)\right\Vert _{+}\leq C_{1}^{N}+C_{2}^{N}\left\Vert
y\right\Vert
\]
follows easily.

Now let $T>\overline{t},y_{0}\in B$ and $\varphi^{\ast}(\cdot)$ be the
corresponding $(\varepsilon,T)$-pseudo-trajectory from Corollary~\ref{cor:2}.
Let us set
\[
y(t)=\left\{
\begin{array}
[c]{l}%
\varphi^{\ast}(t),\ \mbox{if}\ t\not \in \bigcup\limits_{j=0}^{\infty}%
(s_{j}-\overline{t},s_{j}),\\
y_{\varepsilon_{j}}(t-s_{j}+\overline{t},\varphi_{j}(s_{j}-\overline
{t})),\ \mbox{if}\ t\in(s_{j}-\bar{t},s_{j})
\end{array}
\right.
\]
Then $y(\cdot)$ is absolutely continuous function which satisfies
(\ref{vv_eq_14}) and $B\subset\omega(\varphi^{\ast})\subset\omega(y)$. The
theorem is proved.
\end{proof}

\bigskip

We observe that we cannot expect that $F^{N}$ could be replaced by $F$ in
Theorem~\ref{th:2}. The reason is that there are examples in the literature of
connected, compact, chain-recurrent, quasi-invariant sets $B$ which cannot be
included in the $\omega$-limit set of any trajectory of a given differential
equation. Indeed, this fact was shown in \cite{MishST95}. More precisely, let
us consider the Duffing equation
\begin{equation}
\left\{
\begin{array}
[c]{l}%
\dot{x}=y,\\
\dot{y}=x-x^{3}%
\end{array}
\right.  \label{scalareq}%
\end{equation}
The function $V(x,y)=\frac{y^{2}}{2}-\frac{x^{2}}{2}+\frac{x^{4}}{4}$ is
constant along trajectories of (\ref{scalareq}) and any connected component
of
\[
V^{-1}([a,b]),\ \ -\frac{1}{4}\leq a<b\leq\infty,
\]
is compact, invariant, chain recurrent set for the semigroup generated by
(\ref{scalareq}). In particular, it is true for
\[
B=V^{-1}([-\frac{1}{4},0])\cap\{x\geq0\},
\]
and there is no trajectory $\varphi$ of (\ref{scalareq}) such that
$B\subset\omega(\varphi)$.

\bigskip

We give a simple sufficient condition for chain-recurrence.

\begin{lemma}
\label{lem:1} Let $G$ be a strict m-semiflow and for arbitrary $x\in B$ there
exists $t_{x}>0$ such that $x\in G(t_{x},x)$. Then the set $B$ is chain-recurrent.
\end{lemma}

\begin{proof}
Inclusion $x\in G(t_{x},x)$ implies that for all $n\geq1$,%
\[
x\in G(nt_{x},x).
\]
Then for any $t>0$ there exists $n\geq1$ such that $nt_{x}>t$ and $x\in
G(nt_{x},x)$. Thus, $\{x_{1}=x,x_{2}=x\}$ with $t_{1}=nt_{x}$ is the required
$(\varepsilon,t)$-chain. The lemma is proved.
\end{proof}

\bigskip

Further, we will apply Theorem \ref{ChainIncl} to the following partial
differential inclusion:%
\begin{equation}
\left\{
\begin{array}
[c]{c}%
\dfrac{\partial y}{\partial t}-\Delta y\in f(y),\\
y\mid_{\partial\Omega}=0,\\
y(x,0)=y_{0}(x),\ x\in\Omega,
\end{array}
\right.  \label{RDINcl}%
\end{equation}
where $\Omega\subset\mathbb{R}^{n}$ is an open, bounded subset with smooth
boundary $\partial\Omega$ and the multivalued map $f:\mathbb{R}\rightarrow
\mathcal{P}(\mathbb{R})$ satisfies the following assumptions:

\begin{itemize}
\item[$\left(  f1\right)  $] $f$ has non-empty, closed, convex, bounded values.

\item[$\left(  f2\right)  $] $f$ is upper semicontinuous.

\item[$\left(  f3\right)  $] There are $C_{1},C_{2}>0$ such that $\sup_{z\in
f\left(  y\right)  }\left\vert z\right\vert \leq C_{1}+C_{2}\left\vert
y\right\vert .$

\item[$\left(  f4\right)  $] There exist $C_{3},\varepsilon>0$ for which%
\[
zy\leq\left(  \lambda_{1}-\varepsilon\right)  y^{2}+C_{3},\text{ for any }z\in
f\left(  y\right)  ,
\]
being $\lambda_{1}>0$ the first eigenvalue of the operator $-\Delta$ in
$H_{0}^{1}\left(  \Omega\right)  .$
\end{itemize}

Let $V=H_{0}^{1}\left(  \Omega\right)  $, $V^{\prime}=H^{-1}\left(
\Omega\right)  ,\ H=L^{2}\left(  \Omega\right)  $. The operator
$A:V\rightarrow V^{\prime},$ defined by
\[
\langle Au,v\rangle=\int_{\Omega}\nabla u\text{\textperiodcentered}\nabla
vdx,
\]
is continuous and self-adjoint. Also, the operator $-\Delta:D(-\Delta)\subset
H\rightarrow H$ is the subdifferential of the proper, convex lower
semicontinuous function
\[
\psi\left(  u\right)  =\left\{
\begin{array}
[c]{c}%
\frac{1}{2}\int_{\Omega}\left\vert \nabla u\right\vert ^{2}dx\text{, if }u\in
V,\\
+\infty\text{, otherwise,}%
\end{array}
\right.
\]
and the level set given in (\ref{Level}) are compact in $H$. Moreover, it is
easy to see from $\left(  f4\right)  $ that condition (\ref{Diss}) is satisfied.

We define the Nemitski operator $F:H\rightarrow\mathcal{P}(H)$ given by%
\[
F(y)=\{\xi\in H\mid\xi\left(  x\right)  \in f\left(  y\left(  x\right)
\right)  \text{ for a.a. }x\in\Omega\}.
\]
In view of Lemma 6.28 in \cite{KMVY} $F$ satisfies (\ref{vv_eq_3}).

The map $f$ can be written as follows:%
\begin{equation}
f(y)=[\underline{f}(y),\bar{f}(y)],\ \forall y\in\mathbb{R}\text{,}
\label{fCharact}%
\end{equation}
where $\underline{f}:\mathbb{R}\rightarrow\mathbb{R}$ is a lower
semicontinuous single-valued function and $\bar{f}:\mathbb{R\rightarrow R}$ is
an upper semicontinuous single-valued function. Indeed, due to $(f1)$ the map
$f$ can be written as (\ref{fCharact}) with single-valued functions
$\underline{f}$, $\bar{f}$. Since $f$ is upper semicontinuous by $\left(
f2\right)  $, for any $\epsilon>0$ and $y_{0}\in\mathbb{R}$ there exists
$\delta>0$ such that
\[
f(y)\subset O_{\epsilon}(f(y_{0})),\ \text{for all }y\in O_{\delta}(y_{0}).
\]
Therefore, for any $y\in O_{\delta}(y_{0})$ one has%
\[
\underline{f}(y)>\underline{f}(y_{0})-\epsilon,\ \ \bar{f}(y)<\bar{f}%
(y_{0})+\epsilon
\]
and we obtain the required semicontinuity properties.

We recall that the multivalued map $h:\Omega\rightarrow\mathcal{P}%
(\mathbb{R})$ is called measurable if for any open set $O\subset\mathbb{R}$
the inverse
\[
h^{-1}\left(  O\right)  =\{z\in\Omega:h(z)\cap O\not =\varnothing\}
\]
is measurable.

\begin{lemma}
\label{MeasurableComp}If $u:\Omega\rightarrow\mathbb{R}$ is measurable, then
the multivalued map $h:\Omega\rightarrow\mathcal{P}(\mathbb{R})$ given by the
composition $h\left(  x\right)  =f(u(x))$ is measurable.
\end{lemma}

\begin{proof}
Since it is known that any open set of $\mathbb{R}$ is the union of a sequence
of intervals $\left(  a_{n},b_{n}\right)  $, where $a_{n},b_{n}\in
\overline{\mathbb{R}}=\mathbb{R}\cup\{\pm\infty\}$, and $h^{-1}(A\cup
B)=h^{-1}(A)\cup h^{-1}(B)$, it is enough to show that the inverse of any set
$\left(  a,b\right)  $ is Lebesgue measurable. This follows from the
equalities%
\begin{align*}
h^{-1}(\left(  a,b\right)  )  &  =\{x\in\Omega\mid f(u(x))\cap\left(
a,b\right)  \not =\varnothing\}\\
&  =\{x\in\Omega\mid\exists\xi\in f(u(x))\text{ such that }a<\xi<b\}\\
&  =\{x\in\Omega\mid\underline{f}(u(x))<b,\ \bar{f}(u(x))>a\}\\
&  =\{x\in\Omega\mid\underline{f}(u(x))<b\}\cap\{x\in\Omega\mid\ \bar
{f}(u(x))>a\}
\end{align*}
and the fact that the single-valued maps \underline{$g$}$(x)=\underline{f}%
(u(x)),\ \overline{g}\left(  x\right)  =\bar{f}(u(x))$ are measurable.
\end{proof}

\begin{lemma}
\label{ApprSeq}There exists an approximative sequence $F_{N}:H\rightarrow
\mathcal{P}(H)$ having closed, convex, bounded values and satisfying
(\ref{vv_eq_8})-(\ref{vv_eq_10}).
\end{lemma}

\begin{proof}
Taking into account (\ref{fCharact}), we define the Moreau-Yosida
regularization:
\[
f_{N}(x)=\sup_{y\in\mathbb{R}}(\bar{f}(y)-\frac{N}{2}|y-x|^{2}).
\]
Let us prove that this sequence of functions possesses the following properties:

\begin{enumerate}
\item $f_{N}(x)<+\infty,$ for all $N\geq1,\ x\in\mathbb{R}$,

\item for all $x\in\mathbb{R}$ there exist$\ x_{N}$ such that $f_{N}%
(x)=\bar{f}(x_{N})-\frac{N}{2}|x_{N}-x|^{2}$,

\item $\bar{f}(x)\leq f_{N}(x)\leq\bar{f}(x_{N}),\ $for all $N\geq
1,\ x\in\mathbb{R}$,

\item ${x_{N}}\rightarrow x$ $,\bar{f}(x_{N})\rightarrow\bar{f}(x),\ \ f_{N}%
(x)\rightarrow\bar{f}(x),\ \ $as $N\rightarrow\infty,$

\item $|f_{N}(x)|\leq D(1+|x|),\ $\ for any $x\in\mathbb{R}$,

\item $|f_{N}(x)-f_{N}(y)|\leq DN(1+|x|+|y|)|x-y|$, for any $x,y\in\mathbb{R}$
and $N\geq1,$
\end{enumerate}

\noindent where $D>0$ does not depend on $N,x,y$. Indeed, points 1,2 are a
consequence of the Weierstrass theorem and the sublinear growth of $\bar{f}$.
Point 3 is obvious. From the sublinear growth of $\bar{f}$ and the inequality
\[
\bar{f}(x_{N})-\frac{N}{2}|x_{N}-x|^{2}\geq\bar{f}(x)
\]
we deduce that
\begin{equation}
|x_{N}|\leq\bar{C}(1+|x|),\ \label{Boundednessxn}%
\end{equation}
where$\ \bar{C}\ $does\ not depend on $x,N$. These inequalities imply the
convergence $x_{N}\rightarrow x$, and the first inequality and the upper
semicontinuity of $\bar{f}$ imply the convergence $\bar{f}(x_{N}%
)\rightarrow\bar{f}(x)$. Hence, using point 3 we get that $f_{N}%
(x)\rightarrow\bar{f}(x)$. The sublinear growth of $\bar{f}$,
(\ref{Boundednessxn}) and point 3 imply point 5. For proving point 6 we will
use that for any $x,y\in\mathbb{R}$ and their corresponding sequences
$x_{N},y_{N}$ one has%
\[
f_{N}(x)=\bar{f}(x_{N})-\frac{N}{2}|x_{N}-x|^{2}\geq\bar{f}(y_{N})-\frac{N}%
{2}|y_{N}-x|^{2},
\]%
\[
f_{N}(y)=\bar{f}(y_{N})-\frac{N}{2}|y_{N}-y|^{2}\geq\bar{f}(x_{N})-\frac{N}%
{2}|x_{N}-y|^{2}.
\]
Combining these inequalities we obtain
\[
f_{N}(x)-f_{N}(y)=\bar{f}(x_{N})-\bar{f}(y_{N})-\frac{N}{2}|x_{N}-x|^{2}%
+\frac{N}{2}|y_{N}-y|^{2}\leq\frac{N}{2}|x_{N}-y|^{2}-\frac{N}{2}|x_{N}%
-x|^{2},
\]%
\[
f_{N}(y)-f_{N}(x)=\bar{f}(y_{N})-\bar{f}(x_{N})+\frac{N}{2}|x_{N}-x|^{2}%
-\frac{N}{2}|y_{N}-y|^{2}\leq\frac{N}{2}|y_{N}-x|^{2}-\frac{N}{2}|y_{N}%
-y|^{2}.
\]
From the last inequalities and (\ref{Boundednessxn}) we deduce point 6.

For $N\geq1$ let $D_{N+1}$ be the Lipschitz constant of $f_{N}$ in
$[-N-1,N+1]$. Then%
\begin{align*}
f_{N}(x)  &  \leq f_{N}(N)+D_{N+1}\left(  x-N\right)  \text{ if }x\in\lbrack
N,N+1],\\
f_{N}(x)  &  \leq f_{N}(-N)-D_{N+1}(x+N)\text{ if }x\in\lbrack-N-1,-N].
\end{align*}
On the other hand, we know that $f_{N}(x)\leq D+D\left\vert x\right\vert $ for
any $x$. We choose $K_{N}^{+}\geq D_{N+1},\ K_{N}^{-}\geq D_{N+1}\ $ such that
the point of intersection $x_{N}^{+}$ of $f_{N}(N)+K_{N}^{+}\left(
x-N\right)  $ with $D+Dx$ (respectively, $x_{N}^{-}$ of $f_{N}(-N)-K_{N}%
^{-}(x+N)$ with $D-Dx$) belongs to $[N,N+1]$ (respectively, to $[-N-1,-N]$).
We put
\[
{f}^{(N)}(x)=\left\{
\begin{array}
[c]{l}%
D-Dx,\text{ if }x\leq x_{N}^{-},\\
f_{N}(-N)-K_{N}^{-}(x+N)\text{, if }x_{N}^{-}\leq x\leq-N,\\
f_{N}(x),\ \text{if}\ |x|\leq N,\\
f_{N}(N)+K_{N}^{+}(x-N),\ \text{if}\ N\leq x\leq x_{N}^{+},\\
D\left(  1+x\right)  ,\ \text{if}\ x\geq x_{N}^{+}.
\end{array}
\right.
\]
Then:

\begin{enumerate}
\item $\bar{f}(x)\leq{f}^{(N)}(x),\ $for any $N\geq1,\ x\in\mathbb{R},$

\item $|{f}^{(N)}(x)|\leq{D}(1+|x|)$, for any $N\geq1,\ x\in\mathbb{R}$,

\item $\ {f}^{(N)}(x)\rightarrow\bar{f}(x),\ \ $as $N\rightarrow\infty$, for
any $x\in\mathbb{R}$,

\item for any $N\geq1\ $there exists $C(N)>0\ $such that $\left\vert
f^{(N)}(x)-f^{(N)}(y)\right\vert \leq C(N)\left\vert x-y\right\vert $, for all
$x,y\in\mathbb{R}$.
\end{enumerate}

In the same way, for the function $\underline{f}\left(  x\right)  $ we define
a sequence of functions $g^{(N)}\left(  x\right)  $ satisfying:

\begin{enumerate}
\item $\underline{f}\left(  x\right)  \geq g^{(N)}(x),\ $for any $N\geq1$,
$x\in\mathbb{R}$,

\item $|{g}^{(N)}(x)|\leq{D}(1+|x|)$, for any $N\geq1,\ x\in\mathbb{R}$,

\item $\ {g}^{(N)}(x)\rightarrow\underline{f}\left(  x\right)  ,\ \ $as
$N\rightarrow\infty$, for any $x\in\mathbb{R}$,

\item for any $N\geq1\ $there exists $C(N)>0\ $such that $\left\vert
g^{(N)}(x)-g^{(N)}(y)\right\vert \leq C(N)\left\vert x-y\right\vert $, for all
$x,y\in\mathbb{R}$.
\end{enumerate}

We define now the maps $F_{N}:H\rightarrow\mathcal{P}(H)$ by
\[
F_{N}(y)=\{\xi\in H:\xi\left(  x\right)  \in\lbrack{g}^{(N)}\left(
y(x)\right)  ,{f}^{(N)}(y(x))]\text{, for a.a. }x\in\Omega\}.
\]
It follows from \cite[Lemmas 11, 12]{MelVal98} that $F_{N}$ has non-empty,
closed, convex, bounded values and that there exists $c_{N}>0$ such that
(\ref{vv_eq_10}) holds true. Also, it is obvious that (\ref{vv_eq_8}) is
satisfied. Finally, let us prove (\ref{vv_eq_9}). Assume the opposite, that
is, that there exists $y\in H$, $\varepsilon>0$ and a sequence $\xi^{N}\in
F_{N}(y)$ such that%
\begin{equation}
dist(\xi^{N},F(y))>\varepsilon. \label{EpsContr}%
\end{equation}
The multivalued function $x\mapsto f(y(x))$ is measurable by Lemma
\ref{MeasurableComp}. As the functions $\underline{f}\left(  y(x)\right)  ,$
$\bar{f}(y(x)),$ ${g}^{(N)}\left(  y(x)\right)  ,$ ${f}^{(N)}(y(x))$ are
measurable, the map%
\[
\rho_{N}(x)=\left\vert {f}^{(N)}(y(x))-\bar{f}(y(x))\right\vert +\left\vert
{g}^{(N)}\left(  y(x)\right)  -\underline{f}\left(  y(x)\right)  \right\vert
+\frac{1}{N}%
\]
is measurable as well. Let $B(a,c)=[a-c,a+c]$. We define the multivalued map
$\mathcal{P}_{N}:\Omega\rightarrow\mathcal{P}(\mathbb{R})$ by $\mathcal{P}%
_{N}(x)=B(\xi^{N}(x),\rho_{N}(x)).$ By \cite[p.316]{AubinFrank} this map is
measurable. Then the multivalued map $\mathcal{D}_{N}:\Omega\rightarrow
\mathcal{P}(\mathbb{R})$ given \ by $\mathcal{D}_{N}(x)=\mathcal{P}_{N}(x)\cap
f(y(x))$ has non-empty values and is measurable as the intersection of
measurable maps \cite[p.312]{AubinFrank}. Therefore, there exists a measurable
selection $z^{N}(x)\in\mathcal{D}_{N}(x)$ \cite[p.308]{AubinFrank}. Since
$\ {g}^{(N)}(x)\rightarrow\underline{f}\left(  x\right)  ,\ \ {f}%
^{(N)}(x)\rightarrow\bar{f}(x)$, we obtain that $\rho_{N}(x)\rightarrow0$ as
$N\rightarrow\infty$, so%
\[
\left\vert \xi^{N}\left(  x\right)  -z^{N}(x)\right\vert \rightarrow0\text{,
for a.a. }x\in\Omega\text{.}%
\]
The uniform sublinear growth of the functions$\ {g}^{(N)},\ {f}^{(N)}$ implies
that
\[
\left\vert \xi^{N}\left(  x\right)  -z^{N}(x)\right\vert \leq K\left(
1+\left\vert y(x)\right\vert \right)  \text{, for a.a. }x\in\Omega\text{,\ for
any }N\text{,}%
\]
for some constant $K>0$. Hence, Lebesgue's dominated convergence theorem
yields%
\[
v^{N}=\xi^{N}-z^{N}\rightarrow0\text{ in }L^{2}\left(  \Omega\right)  \text{.}%
\]
Hence,%
\[
dist\left(  \xi^{N},F(y)\right)  \leq\left\Vert \xi^{N}-z^{N}\right\Vert
\rightarrow0\text{,}%
\]
which contradicts (\ref{EpsContr}).
\end{proof}

\bigskip

In view of Lemma \ref{ApprSeq} and Corollary \ref{ChainInclCor} we obtain the
following result.

\begin{theorem}
For any solution $\varphi\in\mathcal{K}$ of problem (\ref{RDINcl}) the set
$\omega(\varphi)$ is chain recurrent.
\end{theorem}

\bigskip

\section{Stable and unstable sets of compact quasi-invariant isolated sets
\label{SectionStable}}

In this section we generalize to multivalued semiflows a well known result
about the existence of stable and unstable sets for compact isolated invariant
sets of flows (see \cite{ButlerWaltman}, \cite{HaleWaltman}).

\begin{definition}
A quasi-invariant set $M$ is said to be isolated if there exists a
neighborhood $O$ of $M$ such that $M$ is the maximal quasi-invariant set in
$O$. We call $O$ an isolating neighborhood.
\end{definition}

If the quasi-invariant set $M$ is compact, then in the previous definition we
can replace $O$ by a $\delta$-neighborhood $O_{\delta}(M)$.

We denote the set of all complete trajectories of $\mathcal{K}$ by
$\mathcal{F}$.

\begin{definition}
Let $M$ be an isolated quasi-invariant set. The weakly stable set $W_{w}%
^{s}(M)$ of $M$ is defined by%
\begin{equation}
W_{w}^{s}(M)=\{x\in X\mid\exists\phi\in\mathcal{K}\text{ with }\phi\left(
0\right)  =x\text{ such that }\omega\left(  \phi\right)  \not =\varnothing
,\ \omega\left(  \phi\right)  \subset M\}. \label{DefStableSet}%
\end{equation}
The unstable set $W^{u}(M)$ of $M$ is given by%
\begin{equation}
W^{u}(M)=\{x\in X\mid\exists\varphi\in\mathcal{F}\text{ with }\varphi\left(
0\right)  =x\text{ such that }\alpha\left(  \varphi\right)  \not =%
\varnothing\text{, }\alpha\left(  \varphi\right)  \subset M\}.
\label{DefUnstableSet}%
\end{equation}

\end{definition}

\begin{remark}
The strong stable set $W_{str}^{s}(M)$ of $M$ would be defined by%
\[
W_{s}^{+}(M)=\{x\in X\mid\omega\left(  x\right)  \subset M\},
\]
where
\[
\omega\left(  x\right)  =\{y\in X\mid\exists t_{n}\rightarrow+\infty\text{,
}y_{n}\in G(t_{n},x)\text{ such that }y_{n}\rightarrow y\}.
\]
Of course, in the single-valued case the weak and strong stable sets are the same.
\end{remark}

\begin{theorem}
\label{StableUnstable}Let $\left(  K1\right)  -\left(  K4\right)  $ hold and
let the closure of the positive orbit of $\varphi\in\mathcal{K}$ be compact.
We consider the compact, isolated, quasi-invariant set $M$. If $\omega\left(
\varphi\right)  \cap M\not =\varnothing$ but $\omega\left(  \varphi\right)
\not \subset M$, then
\begin{align}
\omega\left(  \varphi\right)  \cap\{W_{w}^{s}(M)\backslash M\}  &
\not =\varnothing,\label{Stable}\\
\omega\left(  \varphi\right)  \cap\{W^{u}(M)\backslash M\}  &  \not =%
\varnothing. \label{Unstable}%
\end{align}

Moreover, we obtain a function $\phi\in\mathcal{K}$ in the Definition
(\ref{DefStableSet}) such that $\phi\left(  t\right)  \in\omega\left(
\varphi\right)  \cup M$ for all $t\geq0$, and a complete trajectory $\phi$ in
the Definition (\ref{DefUnstableSet}) such that $\phi\left(  t\right)
\in\omega\left(  \varphi\right)  \cup M$ for all $t\in\mathbb{R}$.
\end{theorem}

\begin{proof}
We prove first (\ref{Stable}), that is, the existence of the weakly stable set.

The assumptions $\omega\left(  \varphi\right)  \cap M\not =\varnothing$ but
$\omega\left(  \varphi\right)  \not \subset M$ imply the existence of a
$\delta$-neighborhood $O_{\delta}(M)$ such that the map $\varphi\left(
t\right)  $ enters and leaves it infinitely often as $t\rightarrow+\infty$.
Hence, there exist $0<\tau_{k}<t_{k}$ such that $\tau_{k},t_{k}\rightarrow
+\infty$ and%
\begin{align*}
dist\left(  \varphi\left(  \tau_{k}\right)  ,M\right)   &  =\delta,\\
dist(\varphi\left(  t\right)  ,M)  &  <\delta,\ \text{for all }t\in(\tau
_{k},t_{k}],\\
dist(\varphi\left(  t_{k}),M\right)   &  <\frac{1}{k}.
\end{align*}
We choose $O_{\delta}(M)$ in such a way that $\overline{O_{\delta}(M)}\subset
O_{\varepsilon}(M)$, being $O_{\varepsilon}(M)$ an isolating neighborhood of
$M$.

First, let $t_{k}-\tau_{k}$ be bounded. Then up to \ a subsequence $t_{k}%
-\tau_{k}\rightarrow T>0$. We put $\varphi_{k}(t)=\varphi\left(  t+\tau
_{k}\right)  $. Then, as $\tau_{k}\rightarrow+\infty$, we obtain that%
\begin{align*}
\varphi_{k}(0)  &  =\varphi(\tau_{k})\rightarrow\varphi_{0}\in\omega
(\varphi),\\
dist(\varphi_{0},M)  &  =\delta,
\end{align*}
so $\varphi_{0}\not \in M$. By $\left(  K4\right)  $ there exists
$\overline{\varphi}\in\mathcal{K}$ such that
\[
\varphi_{k}(s_{k})\rightarrow\overline{\varphi}\left(  s\right)  \text{ as
}s_{k}\rightarrow s\text{, }s_{k},s\geq0.
\]
In particular,%
\begin{align*}
\varphi_{k}(0)  &  \rightarrow\overline{\varphi}\left(  0\right)  =\varphi
_{0},\\
\varphi_{k}(t_{k}-\tau_{k})  &  \rightarrow\overline{\varphi}\left(  T\right)
=y\in M.
\end{align*}
Since $M$ is quasi-invariant, there exists $\psi\in\mathcal{K}$ such that
$\psi\left(  0\right)  =y$, $\psi\left(  t\right)  \in M$ for any $t\geq0$.
Using $\left(  K3\right)  $ we can concatenate $\overline{\varphi}$ and $\psi$
by%
\[
\phi\left(  t\right)  =\left\{
\begin{array}
[c]{c}%
\overline{\varphi}(t)\text{ if }0\leq t\leq T,\\
\psi(t-T)\text{ if }t\geq T,
\end{array}
\right.
\]
so that%
\[
\varphi_{0}\in\omega\left(  \varphi\right)  \cap\{W_{w}^{+}(M)\backslash M\}.
\]
It is clear that $\phi\left(  t\right)  \in\omega(\varphi)\cup M$ for any
$t\geq0$.

Second, let $t_{k}-\tau_{k}\rightarrow+\infty$. Since $\tau_{k}\rightarrow
+\infty$, up to a subsequence $\varphi(\tau_{k})\rightarrow y\in\omega
(\varphi)$. We put $\varphi_{k}(t)=\varphi(t+\tau_{k})$. By $\left(
K4\right)  $ there exists $\psi\in\mathcal{K}$ for which $\varphi
_{k}(t)\rightarrow\psi(t)$ uniformly in compact sets of $[0,+\infty)$.
Moreover,%
\begin{align*}
\psi(0)  &  =y\in\omega(\varphi),\\
dist(y,M)  &  =\delta,\\
\psi(t)  &  \in O_{\delta}(M)\text{ for all }t>0.
\end{align*}
The last result is a consequence of $dist(\varphi_{k}\left(  t\right)
,M)<\delta,\ $for all $t\in(0,t_{k}-\tau_{k}].$

Hence, $\omega(\psi)\subset\overline{O_{\delta}(M)}\subset O_{\varepsilon}%
(M)$. As $\psi\left(  t\right)  \in\omega(\varphi)$ for all $t\geq0$,
$\omega\left(  \psi\right)  $ is non-empty and quasi-invariant, so
$\omega\left(  \psi\right)  \subset M$. Thus,
\[
y\in\omega\left(  \varphi\right)  \cap\{W_{w}^{s}(M)\backslash M\}.
\]

Second, we prove (\ref{Unstable}), that is, the existence of the unstable set.

There exist $0<t_{k}<\tau_{k}$ such that $t_{k},\tau_{k}\rightarrow+\infty$
and
\begin{align*}
dist\left(  \varphi\left(  \tau_{k}\right)  ,M\right)   &  =\delta,\\
dist(\varphi\left(  t\right)  ,M)  &  <\delta,\ \text{for all }t\in\lbrack
t_{k},\tau_{k}),\\
dist(\varphi\left(  t_{k}),M\right)   &  <\frac{1}{k}.
\end{align*}

As before, first let $\tau_{k}-t_{k}$ be bounded, so passing to a subsequence
$\tau_{k}-t_{k}\rightarrow T$. We define $\varphi_{k}(t)=\varphi\left(
t+t_{k}\right)  $. We know that
\[
\varphi_{k}(0)=\varphi\left(  t_{k}\right)  \rightarrow x\in M
\]
and by $\left(  K4\right)  $ there exists $\overline{\varphi}\in\mathcal{K}$
such that $\overline{\varphi}(0)=x$ and
\[
\varphi_{k}(s_{k})\rightarrow\overline{\varphi}(s)\text{ if }s_{k}\rightarrow
s\text{, }s_{k},s\geq0\text{.}%
\]
Hence,
\begin{align*}
\varphi_{k}(\tau_{k}-t_{k})  &  =\varphi(\tau_{k})\rightarrow\overline
{\varphi}(T)=\varphi_{0}\in\omega(\varphi),\\
dist\left(  \varphi_{0},M\right)   &  =\delta,
\end{align*}
so $\varphi_{0}\not \in M$ and $T>0$. As $M$ is quasi-invariant, there is a
complete trajectory $\psi$ such that
\begin{align*}
\psi\left(  0\right)   &  =x,\\
\psi(t)  &  \in M\text{ for all }t\leq0.
\end{align*}
Using $\left(  K3\right)  $ we can concatanate $\psi$ and $\overline{\varphi}
$ by%
\[
\phi\left(  t\right)  =\left\{
\begin{array}
[c]{c}%
\psi\left(  t+T\right)  \text{ if }t\leq-T,\\
\overline{\varphi}\left(  t+T\right)  \text{ if }t\geq-T,
\end{array}
\right.
\]
and obtain that $\varphi_{0}\in\omega\left(  \varphi\right)  \cap
\{W^{-}(M)\backslash M\}.$ It is obvious that $\phi\left(  t\right)  \in
\omega\left(  \varphi\right)  \cup M$ for all $t\in\mathbb{R}.$

Assume now that $\tau_{k}-t_{k}\rightarrow+\infty$. As $\tau_{k}%
\rightarrow+\infty$, up to a subsequence $y_{k}=\varphi\left(  \tau
_{k}\right)  \rightarrow y\in\omega\left(  \varphi\right)  $. Let $\varphi
_{k}(t)=\varphi\left(  t+\tau_{k}\right)  $. From \cite[Lemma 13]{CostaValero}
there exists a complete trajectory $\psi$ satisfying $\varphi_{k}%
(t)\rightarrow\psi(t)$ uniformly in bounded sets. Moreover,%
\begin{align*}
\psi(0)  &  =y\in\omega(\varphi),\\
dist(y,M)  &  =\delta,\\
\psi(t)  &  \in O_{\delta}(M)\text{ for all }t<0.
\end{align*}
The last inclusion is a consequence of the fact that $\varphi_{k}(t)\in
O_{\delta}(M)$ for all $t\in\lbrack-\tau_{k}+t_{k},0)$. Thus, $\alpha\left(
\psi\right)  \subset\overline{O_{\delta}(M)}\subset O_{\varepsilon}(M)$. But
$\psi\left(  t\right)  \in\omega(\varphi),$ for all $t,$ implies that
$\alpha\left(  \psi\right)  $ is non-empty and quasi-invariant, so
$\alpha\left(  \psi\right)  \subset M$. Therefore, we have
\[
y\in\omega\left(  \varphi\right)  \cap\{W^{u}(M)\backslash M\}.
\]

\end{proof}

\bigskip

With a similar proof the same result is obtained for the $\alpha$-limit set of
a complete trajectory.

\begin{theorem}
Let $\left(  K1\right)  -\left(  K4\right)  $ hold and let the closure of the
negative orbit of a complete trajectory $\phi\in\mathcal{F}$ be compact. We
consider the compact, isolated, quasi-invariant set $M$. If $\alpha\left(
\phi\right)  \cap M\not =\varnothing$ but $\alpha\left(  \phi\right)
\not \subset M$, then
\begin{align*}
\alpha\left(  \phi\right)  \cap\{W_{w}^{s}(M)\backslash M\}  &  \not =%
\varnothing,\\
\alpha\left(  \phi\right)  \cap\{W^{u}(M)\backslash M\}  &  \not =\varnothing.
\end{align*}

\end{theorem}

The existence of stable and unstable sets given in Theorem
\ref{StableUnstable} allows us to prove the existence of homoclinic structures
inside $\omega$-limit sets of trajectories. This is the opposite situation to
that considered in \cite{a-cccl} (see also \cite{CostaValero}), in which Morse
decompositions, characterized by the absense of homoclinic trajectories, are constructed.

We say that there exists a connection from the set $M$ to the set $N$ if there
are $x\not \in M\cup N$ and a bounded complete trajectory $\phi\in\mathcal{F}$
satisfying $\phi\left(  0\right)  =x$ and $\omega\left(  \phi\right)  \subset
N$, $\alpha\left(  \phi\right)  \subset M$, being $\omega\left(  \phi\right)
$, $\alpha\left(  \phi\right)  $ non-empty. When $M$ and $N$ are disjoint the
connection is called heteroclinic, whereas if $M=N$, it is called homoclinic.
It is also said that $M$ is chained to $N$.

A finite number of pairwise disjoint sets $\{M_{1},...,M_{k}\}$, $k\geq1$, is
said to be cyclically chained to each other (or a cyclical chain) if for any
$j\in\{1,...,k\}$ there exists a connection from $M_{j}$ to $M_{j+1}$, where
$M_{k+1}=M_{1}$. When $k=1$, $M_{1}$ is just chained to itself. A cyclical
chain is called also an homoclinic structure.

Let us consider the $\omega$-limit set $\omega\left(  \varphi\right)  $ of a
trajectory $\varphi\in\mathcal{K}$ having a compact positive orbit. Since this
set is quasi-invariant, there exists a set of bounded complete trajectories
$\mathbb{K}_{\varphi}$ such that%
\[
\omega\left(  \varphi\right)  =\{\phi\left(  0\right)  \mid\phi\in
\mathbb{K}_{\varphi}\}.
\]
We define the set $\Omega_{\varphi}$ by%
\[
\Omega_{\varphi}=\cup_{\phi\in\mathbb{K}_{\varphi}}\omega(\phi).
\]
The following lemma extends to the multivalued case Proposition 3.3 in
\cite{Thieme}.

\begin{lemma}
\label{Cyclic}Let $\left(  K1\right)  -\left(  K4\right)  $ hold and let the
closure of the positive orbit of $\varphi\in\mathcal{K}$ be compact. Let
$\Omega_{\varphi}\subset M=\cup_{j=1}^{m}M_{j}$, $m\geq1$, where $M_{j}%
\subset\omega\left(  \varphi\right)  $ are compact, pairwise disjoint,
isolated, quasi-invariant subsets.

Then either $\omega\left(  \varphi\right)  =M_{1}$ (so $m=1$) or there exists
a cyclical chain $\{\widetilde{M}_{1},...,\widetilde{M}_{k}\}$, $1\leq k\leq
m$, where each $\widetilde{M}_{i}$ is equal to some $M_{j}$. Moreover, the
connections in this chain belong entirely to $\omega\left(  \varphi\right)  $.
\end{lemma}

\begin{proof}
If $\omega\left(  \varphi\right)  \not =M_{1}$, we choose some $M_{i}$ and
name it $\widetilde{M}_{1}$. Since $\omega\left(  \phi\right)  \backslash
\widetilde{M}_{1}\not =\varnothing$, Theorem \ref{StableUnstable} and
$M_{i}\subset\omega\left(  \varphi\right)  $ imply the existence of a complete
trajectory $\phi$ such that $\phi\left(  \mathbb{R}\right)  \subset
\omega\left(  \varphi\right)  $, $\phi\left(  0\right)  \not \in
\widetilde{M}_{1}$ and $dist\left(  \phi\left(  t\right)  ,\widetilde{M}%
_{1}\right)  \rightarrow0$ as $t\rightarrow-\infty$. By Lemma \ref{PropOmega}
the set $\omega\left(  \phi\right)  $ is connected, so it belong to one of the
sets $M_{j}$, renamed $\widetilde{M}_{2}$. If $\widetilde{M}_{2}%
=\widetilde{M}_{1}$, we are done. Otherwise, we have obtained a connection
from $\widetilde{M}_{1}$ to $\widetilde{M}_{2}$. Further, arguing in the same
way we get a connection from $\widetilde{M}_{2}$ to some $\widetilde{M}_{3}$.
If either $\widetilde{M}_{3}=\widetilde{M}_{1}$ or $\widetilde{M}%
_{3}=\widetilde{M}_{2}$, we have finished. In other case we continue in the
same way until in a finite number of steps we obtain the desired chain.
\end{proof}

\bigskip

\section{Convergence to equilibria for reaction-diffusion equations without
uniqueness}

In a bounded domain $\Omega\subset\mathbb{R}^{n}$, $1\leq n\leq3,$ with
sufficiently smooth boundary $\partial\Omega$ we consider the problem
\begin{equation}
\left\{
\begin{array}
[c]{l}%
u_{t}-\Delta u+f(u)=h,\quad x\in\Omega,\ t>0,\\
u|_{\partial\Omega}=0,\ t>0,\\
u\left(  0,x\right)  =u_{0}\left(  x\right)  \text{, }x\in\Omega,
\end{array}
\right.  \label{RD}%
\end{equation}
where
\begin{equation}%
\begin{array}
[c]{c}%
h\in L^{\infty}(\Omega),\ \\
f\in C(\mathbb{R}),\\
|f(u)|\leq\alpha(1+|u|^{p-1}),\quad\forall u\in\mathbb{R},\\
f(u)u\geq\beta u^{p}+\gamma,\ \forall u\in\mathbb{R},
\end{array}
\label{CondRD}%
\end{equation}
with $2\leq p\leq4$, $\alpha,\beta,\gamma>0$.

Let $H=L^{2}\left(  \Omega\right)  $, $V=H_{0}^{1}(\Omega)$, whereas
$\Vert\cdot\Vert,\ (\cdot,\cdot)$ will be the norm and the scalar product in
$L^{2}(\Omega)$.

A function $u\in L_{loc}^{2}(0,+\infty;V)\bigcap L_{loc}^{p}(0,+\infty
;L^{p}(\Omega))$ is called a weak solution of (\ref{RD}) on $(0,+\infty)$ if
for all $T>0\,,\ v\in V,\,\eta\in C_{0}^{\infty}(0,T),$
\[
-\int\limits_{0}^{T}(u,v)\eta_{t}dt+\int\limits_{0}^{T}\left(  (u,v)_{V}%
+(f(u),v)-(h,v)\right)  \eta dt=0.
\]
A weak solution is called a strong one if, moreover, $u\in L^{\infty}\left(
0,T;V\right)  ,\ \dfrac{du}{dt}\in L^{2}\left(  0,T;H\right)  ,\ $for any
$T>0$. Any strong solution $u$ satisfies $u\in L^{2}\left(  0,T;D\left(
A\right)  \right)  \cap C([0,+\infty);V).$

It is well known \cite[p.284]{ChepVishikBook} that for any $u_{0}\in H$ there
exists at least one weak solution of (\ref{RD}) with $u(0)=u_{0}$ (which might
be non unique) and that any weak solution of (\ref{RD}) belongs to $C\left(
[0,+\infty);H\right)  $.

The aim of this section is to prove, using Lemma \ref{Cyclic}, that the
$\omega$-limit set of every weak solution is a stationary point. For this
purpose we need to recall some results about the properties of weak solutions
and the structure of the global attractor for (\ref{RD}) given in
\cite{KapKasVal15}.

Let $\mathcal{K\subset}C(\mathbb{R}^{+},H)$ be the set of all weak solutions
of problem (\ref{RD}) with initial condition in $H$. This set satisfies
properties $\left(  K1\right)  -\left(  K4\right)  $, so it generates the
strict m-semiflow $G:\mathbb{R}^{+}\times H\rightarrow\mathcal{P}(H)$ by
(\ref{141014-1640}). Moreover, $G$ possesses the global compact invariant
attractor $\Theta$ (see the definition in Section \ref{Incl}). Therefore,
every positive orbit $\gamma^{+}(\varphi)$, $\varphi\in\mathcal{K},$ has
compact closure in $H$ and by Lemma \ref{PropOmega}\ its $\omega$-limit set
$\omega(\varphi)$ is nonempty, compact, quasi-invariant, connected and
$\omega(\varphi)\subset\Theta$.

Further, we give an insight into the structure of the global attractor in
terms of bounded complete trajectories. A complete trajectory $\phi$ is said
to be bounded if $\cup_{t\in\mathbb{R}}\phi\left(  t\right)  $ is a bounded
set. We denote by $\mathbb{K}$ the set of all bounded complete trajectories.
Then the global attractor is characterized by%
\[
\Theta=\left\{  \phi(0)\,:\,\phi(\cdot)\in\mathbb{K}\right\}  .
\]
We can give a more detail description of $\Theta$ in terms of the unstable and
weakly stable sets of the stationary points. We denote by $\mathfrak{R}$ the
set of all stationary points of (\ref{RD}), i.e., the points $u\in V$ such
that%
\begin{equation}
-\Delta u+f(u)=h\text{ in }H^{-1}\left(  \Omega\right)  . \label{Stationary}%
\end{equation}
It is known \cite[Lemmas 12, 14 and Theorem 13]{KKV} that $\mathfrak{R}%
\not =\varnothing$ and that the following properties are equivalent:

\begin{enumerate}
\item $u_{0}\in\mathfrak{R;}$

\item $u_{0}\in G\left(  t,u_{0}\right)  $ for all $t\geq0;$

\item The function $u\left(  t\right)  \equiv u_{0}$ belongs to $\mathcal{K}.$
\end{enumerate}

We define now the sets:%
\begin{equation}%
\begin{array}
[c]{c}%
M^{s}(\mathfrak{R})=\left\{  z\,:\,\exists\phi(\cdot)\in\mathbb{K}%
,\,\ \phi(0)=z,\,\,dist(\phi(t),\mathfrak{R})\rightarrow0,\,\ t\rightarrow
+\infty\right\}  ,\\
M^{u}(\mathfrak{R})=\left\{  z\,:\,\exists\phi(\cdot)\in\mathcal{F}%
,\,\ \phi(0)=z,\,\ dist(\phi(t),\mathfrak{R})\rightarrow0,\,\ t\rightarrow
-\infty\right\}  .
\end{array}
\label{M}%
\end{equation}
$M^{u}(\mathfrak{R})$ is the unstable set of $\mathfrak{R}$. $M^{s}%
(\mathfrak{R})$ is the weakly stable set of $\mathfrak{R}$ but considering
only bounded complete trajectories. In the definition of $M^{u}(\mathfrak{R})$
we can replace $\mathcal{F}$ by $\mathbb{K}$, because the positive orbit of
every complete trajectory $\phi$ is bounded.

\begin{lemma}
\cite[Theorems 4, 5]{KapKasVal15} The global attractor $\Theta$ is bounded in
$L^{\infty}\left(  \Omega\right)  $, compact in $V$ and%
\[
\Theta=M^{u}(\mathfrak{R})=M^{s}(\mathfrak{R}).
\]
Moreover, every weak solution $u\left(  \text{\textperiodcentered}\right)  $
with $u\left(  0\right)  \in\Theta$ is a strong solution.
\end{lemma}

In fact, any bounded complete trajectory $\phi$ satisfies the convergences
given in (\ref{M}), that is,%
\begin{align*}
dist(\phi(t),\mathfrak{R})  &  \rightarrow0,\,\ t\rightarrow+\infty,\\
dist(\phi(t),\mathfrak{R})  &  \rightarrow0,\,\ t\rightarrow-\infty.
\end{align*}
This follows from the fact that every solution inside the global attractor is
strong, and then a Lyapunov function exists (see the proof of Theorem 37 in
\cite{KKV}). Hence, $\omega\left(  \phi\right)  \subset\mathfrak{R}$,
$\alpha\left(  \phi\right)  \subset\mathfrak{R}$. In the particular, case when
there exists a finite number of stationary points, as the sets $\omega\left(
\phi\right)  ,\ \alpha\left(  \phi\right)  $ are connected, they have to be
equal to one of the stationary points. In such a case the global attractor
consists of the stationary points and all bounded complete trajectories
connecting them.

We recall that a Lyapunov function $t\mapsto E\left(  u\left(  t\right)
\right)  $ is strictly decreasing if $u\left(  \text{\textperiodcentered
}\right)  $ is not a stationary point. Therefore, if there exists a connection
from a stationary point $e_{1}$ to the stationary point $e_{2}$, then
necessarily $E\left(  e_{2}\right)  <E\left(  e_{1}\right)  .$

We are now in position to prove the main result of this section.

\begin{theorem}
\label{OmegaStationary}Let the number of stationary points be finite. Then any
$\varphi\in\mathcal{K}$ satisfies $\omega\left(  \varphi\right)
=e\in\mathfrak{R}$.
\end{theorem}

\begin{proof}
As we have seen in Section \ref{SectionStable} there exists a set of bounded
complete trajectories $\mathbb{K}_{\varphi}$ such that
\[
\omega\left(  \varphi\right)  =\{\phi\left(  0\right)  \mid\phi\in
\mathbb{K}_{\varphi}\}.
\]
In view of the previous arguments, the set $\Omega_{\varphi}=\cup_{\phi
\in\mathbb{K}_{\varphi}}\omega(\phi)$ belongs to $\mathfrak{R}$. Put then
$M=\Omega_{\varphi}=\cup_{j=1}^{m}e_{j}\subset\mathfrak{R}$ in Lemma
\ref{Cyclic}. Therefore, the sets $M_{j}=e_{j}$, $j=1,...,m$, are the
stationary points in $\omega\left(  \varphi\right)  $. Since there is a finite
number of them, the sets $M_{j}$ are compact, pairwise disjoint and
quasi-invariant. In order to show that they are also isolated, we choose
$\delta>0$ such that $O_{\delta}(M_{j})\cap O_{\delta}(M_{j})=\varnothing$ if
$i\not =j$. Let us assume that there exists a quasi-invariant set $N$ inside
of one $O_{\delta}(e_{j})$ but $N\not =e_{j}$. Therefore, there is a bounded
complete trajectory $\phi$ satisfying $\phi\left(  0\right)  \not =e_{j}$ and
$\phi\left(  \mathbb{R}\right)  \subset N$. We know that $\phi\left(
t\right)  $ has to converge to a stationary point if either $t\rightarrow
+\infty$ or $t\rightarrow-\infty$. Since the only possible point is $e_{j}$,
there is a connection from $e_{j}$ to itself. This would lead to the
contradiction $E\left(  e_{j}\right)  <E(e_{j})$. Thus, the sets $M_{j}$ are isolated.

Finally, if $\omega\left(  \varphi\right)  $ was not equal to a stationary
point, there would exist by Lemma \ref{Cyclic} a cyclic chain in $M$. This is
not possible again by the decreasing property of the Lyapunov function
$E\left(  \phi\left(  t\right)  \right)  $ along the connecting complete
trajectories $\phi$. We deduce that $\omega\left(  \varphi\right)
=e\in\mathfrak{R}$.
\end{proof}

\begin{remark}
The condition $h\in L^{\infty}(\Omega)$ is crucial in this theorem. If $h\in
L^{2}\left(  \Omega\right)  $, then the question about the structure of
$\omega\left(  \varphi\right)  $ for $\varphi\in\mathcal{K}$ remains open.
Nevertheless, if $2\leq p\leq3$, then in \cite{KapKasVal15} it is proved that
every weak solution is regular and for such solutions a Lyapunov function also
exists. Therefore, under this stronger condition on $p$ the result
$\omega\left(  \varphi\right)  =e\in\mathfrak{R}$ is also true.
\end{remark}

We will consider finally a Chafee-Infante problem in which the number of
stationary points is known to be finite:%
\begin{equation}
\left\{
\begin{array}
[c]{c}%
\dfrac{\partial u}{\partial t}-\dfrac{\partial^{2}u}{\partial x^{2}%
}=f(u),\ t>0,\ x\in\left(  0,1\right)  ,\\
u(t,0)=0,\ u(t,1)=0,\\
u(0,x)=u_{0}(x).
\end{array}
\right.  \label{problemainicial}%
\end{equation}
The function $f$ satisfies the conditions in (\ref{CondRD}) and also the
following ones:

\begin{enumerate}
\item $f(0)=0$;

\item $f^{\prime}\left(  0\right)  >0$ exists and is finite;

\item $f$ is strictly concave if $u>0$ and strictly convex if $u<0.$
\end{enumerate}

It was proved in \cite{Caballero} that if $n^{2}\pi^{2}<f^{\prime}\left(
0\right)  \leq\left(  n+1\right)  ^{2}\pi^{2}$, where $n\geq0$, $n\in
\mathbb{Z}$, then there are exactly $2n+1$ stationary points. In this case
$h\equiv0$. Therefore, Theorem \ref{OmegaStationary} implies the following result.

\begin{theorem}
For any weak solution $\varphi\in\mathcal{K}$ of problem
(\ref{problemainicial}) we have that $\omega\left(  \varphi\right)
=e\in\mathfrak{R}$.
\end{theorem}

\bigskip

\textbf{Acknowledgments.}

The first two authors were partially supported by the State Fund for
Fundamental Research of Ukraine. The third author was partially supported by
Spanish Ministry of Economy and Competitiveness and FEDER, projects
MTM2015-63723-P and MTM2016-74921-P.

We would like to thank the anonymous referees for their useful remarks.

\end{document}